\documentclass[a4paper, 12pt, final]{amsart}
\usepackage[OT1]{fontenc}
\usepackage[latin1]{inputenc}
\usepackage{amssymb,latexsym}
\usepackage{type1cm,ae}
\usepackage[notref,notcite]{showkeys}
\usepackage{graphicx}
\usepackage{xspace}
\usepackage{setspace}
\usepackage[english]{babel}

\theoremstyle{plain}		\newtheorem{theorem}{Theorem}
				
				\newtheorem{lemma}{Lemma}
				\newtheorem{thmA}{Theorem}
				
\theoremstyle{definition}	
\theoremstyle{remark}		

\newcommand*{\bR}{\ensuremath{\mathbb{R}}}

\newcommand*{\bC}{\ensuremath{\mathbb{C}}}

\newcommand*{\loc}{\mathrm{loc}}
\newcommand*{\closure}[1]{\overline{#1}}
\newcommand*{\bdary}[1]{\partial #1}

\newcommand*{\bilip}{bi-Lipschitz}

\newcommand*{\Wert}{\mathord{\mbox{|\kern-1.5pt|\kern-1.5pt|}}}
\newcommand*{\ie}{\mbox{i.e.}\xspace}

\newcommand*{\OO}{\mathcal{O}}

\DeclareMathOperator{\dist}{dist}
\DeclareMathOperator{\diam}{diam}

\DeclareMathOperator{\modulus}{mod}
\DeclareMathOperator{\capacity}{cap}

\title[Mappings of finite distortion: Exponential Cusp Domain]{Mappings of finite distortion: formation of Exponential Cusp}

\author{Changyu Guo}
\address[Changyu Guo]{Department of Mathematics and Statistics, University of Jyv\"askyl\"a, P.O. Box 35, FI-40014 University of Jyv\"askyl\"a, Finland}
\email{guocybnu@gmail.com}


\begin{document}
\begin{abstract}
We consider a quasi-convex planar domain $\Omega$ with a rectifiable boundary but containing an  exponential cusp and show that there is no homeomorphism $f\colon\bR^2\to\bR^2$ of finite distortion with $\exp(\lambda K)\in L_{loc}^1(\bR^2)$ for some $\lambda>0$ such that $f(B)=\Omega$. On the other hand, if we only require that $K_f(x)\in L^p_{loc}(\bR^2)$, then such an $f$ exists.
\end{abstract}

\maketitle

\section{Introduction}\label{sec:first}

The concept of quasidisk is central in the theory of planar quasiconformal mappings; see, for example, \cite{a06,aim09,g82,lv73}. One calls a Jordan domain $\Omega\subset\bR^2$ a quasidisk if it is  the image of the unit disk under a quasiconformal mapping $f\colon\bR^2\to\bR^2$ of the entire plane. If $f$ is $K$-quasiconformal, we say that $\Omega$ is a $K$-quasidisk. 

A substantial part of the theory of quasiconformal mapping has recently been shown to extend in a natural form to the setting of mappings of locally exponentially integrable distortion~\cite{agrs10,aim09,iko02,ko03,oz05}. However, very little is known about the analogues of the concept of a quasidisk. For the model domain
\begin{equation}\label{eq:cuspdomain1}
  \Omega_s=\{(x_1,x_2)\in\bR^2:0<x_1<1,|x_2|<x_1^{1+s}\}\cup B(x_s,r_s),
\end{equation}
where $x_s=(s+2,0)$ and $r_s=\sqrt{(s+1)^2+1}, s>0$, the situation is rather well understood: $\Omega_s=f(B(0,1))$ under a homeomorphism with locally $\lambda$-exponentially integrable distortion if $\lambda< 2/s$, but this cannot happen when $\lambda> 2/s$, see~\cite{kt10}. Moreover, if $f$ is additionally required to be quasiconformal in $B(0,1)$, then the critical bound for $\lambda$ is $1/s$, see
~\cite{kt07}. While for the domain 
\begin{equation}\label{eq:cuspdomain2}
  \Delta_s=  B(x'_s,r'_s)\setminus\{(x_1,x_2)\in\bR^2:x_1>0,|x_2|<x_1^{1+s}\}  ,
\end{equation}
where $x'_s=(-s,0)$ and $r'_s=\sqrt{(s+1)^2+1}, s>0$, it is proven in~\cite{gkt10} $\Delta_s=f(B(0,1))$ under a homeomorphism with locally $\lambda$-exponentially integrable distortion for $\lambda< 2/s$, but this cannot happen if $f$ is additionally required to be quasiconformal in $B(0,1)$. One may guess the reason for this is that the domain $\Omega_s$ is quasi-convex while $\Delta_s$ is not. However, this turns out not to be the case.

For the exponential cusp domain  
\begin{equation}\label{eq:Lipdomain}
  \Omega=\{(x_1,x_2)\in\bR^2:0<x_1<1, |x_2|< e^{-\frac{1}{x_1}}\}\cup B(x_0,r_0),
\end{equation}
where $x_0=2, r_0=\sqrt{1+\frac{1}{e^2}}$, as the image of the unit disk $B$ under a planar homeomorphism of finite distortion, we have the following result.

\begin{thmA}\label{thm:thma}
   There is no homeomorphism $f\colon\bR^2\to\bR^2$  of finite distortion with $\exp(\lambda K)\in L_{loc}^1(\bR^2)$ for some $\lambda>0$ such that $f(B)=\Omega$. While for any $p>0$, there is a homeomorphism $f\colon\bR^2\to\bR^2$  of finite distortion with $K_f(x)\in L^p_{loc}(\bR^2)$ such that $f(B)=\Omega$.
\end{thmA}

\section{Notation and Definitions}\label{sec:notdef}
We sometimes associate the plane $\bR^2$ with the complex plane $\bC$ for convenience. The closure of a set $U\subset\bR^2$ is denoted $\closure{U}$ and the boundary $\bdary{U}$. The open disk of radius $r>0$ centered at $x\in\bR^2$ is denoted $B(x,r)$ and in the case of the unit disk we omit the centre and the radius, writing $B:=B(0,1)$.
The symbol~$\Omega$ always refers to a domain, \ie a connected and open subset of~$\bR^2$. We call a homeomorphism $f\colon\Omega\to f(\Omega)\subset\bR^2$ a homeomorphism of finite distortion if $f\in W_{\loc}^{1,1}(\Omega;\bR^2)$ and
\begin{equation}\label{eq:disteq}
  \|Df(x)\|^2\leq K(x)J_f(x) \text{ a.e. in } \Omega,
\end{equation}
for some measurable function~$K(x)\geq1$ that is finite almost everywhere. In the distortion inequality~\eqref{eq:disteq}, $Df(x)$ is the formal differential of~$f$ at the point~$x$ and $J_f(x):=\det Df(x)$ is the Jacobian. The norm of~$Df(x)$ is defined as 
\begin{equation*}
  \|Df(x)\|:=\max_{e\in\bdary{B}} |Df(x)e|.
\end{equation*}
For a homeomorphism of finite distortion it is convenient to write $K_f$ for the optimal distortion function. This is obtained by setting $K_f(x)=\|Df(x)\|^2/J_f(x)$ when $Df(x)$ exists and $J_f(x)>0$, and $K_f(x)=1$ otherwise. The distortion of~$f$ is said to be locally $\lambda$-exponentially integrable if $\exp(\lambda K_f(x))\in L_{\loc}^1(\Omega)$, for some $\lambda>0$. Note that if we assume $K_f(x)$ to be bounded, we recover the class of quasiconformal mappings (cf.~\cite{lv73}); recall that $J_f\in L_{loc}^{1}(\Omega)$ for each homeomorphism $f\in W_{loc}^{1,1}(\Omega;\bR^2)$(cf.~\cite{aikm09}).

Next we define the two central tools for us -- the modulus of a path family and the capacity. Let~$E$ and~$F$ be subsets of $\closure{\Omega}$. We denote by $\Gamma(E,F,\Omega)$ the path family consisting of all locally rectifiable paths joining~$E$ to~$F$ in~$\Omega$. A Borel function $\rho\colon\bR^2\to[0,\infty\mathclose]$ is said to be admissible for $\Gamma(E,F,\Omega)$ if $\int_\gamma\rho\,ds\geq1$ for all $\gamma\in\Gamma(E,F,\Omega)$. The modulus of a path family $\Gamma:=\Gamma(E,F,\Omega)$ is defined as
\begin{multline*}
  \modulus(\Gamma):=\inf\Big\{\int_{\Omega}\rho^2(x)\,dx :  \rho\colon\bR^2\to[0,\infty\mathclose] \text{ is an admissible }\\
  \text{Borel function for } \Gamma \Big\}.
\end{multline*}
By $\modulus_{K_f(x)}(\Gamma)$ we mean the $K_f(x)$-weighted modulus, where instead of $\int\rho^2(x)\,dx$ we take the infimum over $\int\rho^2(x)K_f(x)\,dx$.

Let~$E$ and~$F$ be disjoint compact sets in a domain~$\Omega$. Let~$\omega$ be measurable with $0\leq\omega(x)\leq1$ almost everywhere. The $\omega$-weighted $p$-capacity of the pair $(F,E)$ with respect to~$\Omega$ is defined to be
\begin{multline*}
  \capacity^p_{\omega}(F,E;\Omega):=
  \inf\Big\{\int_{\Omega}|\nabla u(x)|^p\omega(x)\,dx : u\in C(\Omega)\cap W_{\loc}^{1,1}(\Omega), \\
  u\leq0 \text{ on } F \text{ and } u\geq1 \text{ on } E \Big\}.
\end{multline*}

Finally we note that when we write $f(x)\asymp g(x)$, we mean that $g(x)/c\leq f(x)\leq cg(x)$ is satisfied for all~$x$ with some fixed constant $c\geq 1$. The notation $f(x)=g(x)+\OO(|x|^n)$ means that for some fixed constant $C>0$ we have $|f(x)-g(x)|\leq C|x|^n$ when $|x|$ is small.

\section{Auxiliary results}

We begin by introducing the following lemma,  whose proof can be found in ~\cite{kt10}. 

\begin{lemma}\label{lemma:capala1}
  Let $E\subset\closure{B}$ be a continuum such that $E\subset B(x_0,1/6)$ for some $x_0\in\bdary{B}$ and $F:=\closure{B}(0,1/4)$. Suppose that $v\in W^{1,1}(B)$ is continuous and satisfies: $v=0$ on~$F$ and $\lim_{y\to x} v(y)\geq 1$ for every $x\in E$. If $L:=\int_B\exp(\lambda K)<\infty$, for some measurable function $K(x)\geq 1$, then
\begin{equation}
  \int_B\frac{\vert \nabla v\vert^2}{K}\geq C\lambda\bigg(\log\frac{\sqrt{4L/\pi}}{\diam E}\bigg)^{-2}.
\end{equation}
\end{lemma}

\begin{lemma}\label{lemma:lipest1}
   Let $\Omega$ be a domain as defined in~\eqref{eq:Lipdomain} and $F\subset\Omega$ a compact set. Set $d:=\min\{1,\dist(0,F)\}$ and let $0<r<d/2$. If $E\subset\{(x_1,x_2)\in\closure{\Omega}:0\leq |x|\leq r\}$, then there is a Lipschitz function~$u$ on~$\Omega$ such that $u=0$ on~$F$,  $\lim_{y\to x} u(y)=1$ for every $x\in E$ and
\begin{equation}
  \int_{\Omega}\vert\nabla u\vert^2\,dx\leq Cr^s
\end{equation}
for any $s>0$.
\end{lemma}
\begin{proof}
 Define
\begin{equation*}
  u(x_1,x_2)=
  \begin{cases}
     1 & \text{if } x_1\leq r\\
     1-\int_r^{x_1}\frac{dt}{e^{-\frac{1}{t}}}\big(\int_r^{d/2}\frac{dt}{e^{-\frac{1}{t}}}\big)^{-1} & \text{if } r<x_1\leq d/2\\
     0 & \text{if } x_1>d/2,
\end{cases}
\end{equation*}
and compute
\begin{align*}
  \int_{\Omega_s}|\nabla u(x)|^2\,dx&=
\Big(\int_r^{d/2}\frac{dt}{e^{-\frac{1}{t}}}\Big)^{-2}\int_r^{d/2}\int_{-e^{-\frac{1}{t}}}^{e^{-\frac{1}{t}}}\frac{1}{e^{-\frac{2}{x_1}}}dx_2dx_1\\
  &=\Big(\int_r^{d/2}\frac{dt}{e^{-\frac{1}{t}}}\Big)^{-1}\leq Cr^s
\end{align*}
for any $s>0$.
\end{proof}

\begin{lemma}\label{lemma:diamarvio}
  Let $f\colon\bR^2\to\bR^2$ be a homeomorphism of finite distortion such that $\exp(\lambda K)\in L_{\loc}^1(\bR^2)$ for some $\lambda>0$, and $f(B)=\Omega$. Let $E'_t=\{x\in\bdary{\Omega}:|x|\leq t\}$ and $E_t=f^{-1}(E'_t)$. Then for all $\varepsilon>0$ there exists $t_0>0$ such that for some positive constants~$C$ and~$\tilde{C}$
\begin{equation}\label{eq:diamarvio1}
  \diam E_t\geq C\exp\Big(\frac{-\tilde{C}}{(\diam E'_t)^{\frac{1+\varepsilon}{\lambda}}}\Big)
\end{equation}
for all $0<t<t_0$.
\end{lemma}
The proof of this lemma also can be found in~\cite{kt10}.
\section{Main proofs}
\begin{theorem}\label{thm:partulos}
  There is no homeomorphism $f\colon\bR^2\to\bR^2$  of finite distortion with $\exp(\lambda K)\in L_{loc}^1(\bR^2)$ for some $\lambda>0$ such that $f(B)=\Omega$.
\end{theorem}
\begin{proof}
We prove it by contradiction. Suppose such a homeomorphism $f$ exists and let $\varepsilon>0$. Define $E'_t=\{x\in\bdary{\Omega}:|x|\leq t\}$, $F=\closure{B}(0,1/4)$ and set $E_t=f^{-1}(E'_t)$, $F'=f(F)$.
From Lemma~\ref{lemma:diamarvio} we obtain $t_0>0$ such that for some positive constants $C$ and $\tilde{C}$
\begin{equation}\label{eq:epayht1}
  \diam E_t\geq 
  C\exp\Big(\frac{-\tilde{C}}{(\diam E'_t)^{\frac{1+\varepsilon}{\lambda}}}\Big)
\end{equation}
for all $0<t<t_0$. As~$f$ is an homeomorphism we may assume, by making $t_0$ smaller if necessary, that $E_t\subset B(f^{-1}(0),1/6)$ for all $0<t<t_0$ and $\diam E'_{t_0}<\dist(0,F')$.
By applying Lemma~\ref{lemma:lipest1} we obtain a Lipschitz function~$u$ on~$\Omega$ such that $u=1$ on~$E'_t$, $u=0$ on~$F'$, and
\begin{equation}\label{eq:capylaarv1}
  \int_{\Omega}|\nabla u|^2\,dx\leq Ct^s
\end{equation}
for any $s>0$.

Set $v=u\circ f$ and recall that if $f\in W_{\loc}^{1,1}(B,\bR^2)$ is a homeomorphism, then for each non-negative measurable function $w\colon\bR^2\to\bR$,
\begin{equation}\label{eq:mvkaava1}
  \int_B w\circ f\, |J_f|\leq\int_{f(B)}w.
\end{equation}
Now, as~$f$ is a homeomorphism and $f\in W_{\loc}^{1,1}(\bR^2,\bR^2)$ we know that $\nabla v$ exists almost everywhere, $|\nabla v|=|Df||\nabla u\circ f|$ is locally integrable, and thus $v\in C(B)\cap W_{\loc}^{1,1}(B)$. Also, $v=1$ on~$E_t$ and~$v=0$ on~$F$, so~$v$ is clearly admissible for $\capacity_{1/K}(F,E_t;B)$. Since~$f$ has finite exponentially integrable distortion, we can use the distortion inequality together with~\eqref{eq:mvkaava1} to obtain the estimate
\begin{equation*}
  \int_{B}\frac{|\nabla v|^2}{K}\leq\int_{B}|\nabla u\circ f|^2 J_f\leq\int_{\Omega}|\nabla u|^2
\end{equation*}
and thus with~\eqref{eq:capylaarv1} we readily have the inequality
\begin{equation}\label{eq:capylaarv1b}
  \capacity_{1/K}(F,E_t;B)\leq Ct^s.
\end{equation}
Applying Lemma~\ref{lemma:capala1} we obtain
\begin{equation*}
  C_1\lambda\Big(\log\frac{C_2}{\diam E_t}\Big)^{-2}
  \leq\capacity_{1/K}(F,E_t;B)
\end{equation*}
and combining with~\eqref{eq:capylaarv1b} gives the estimate
\begin{equation}\label{eq:valepayht}
   C_1\Big(\log\frac{C_2}{\diam E_t}\Big)^{-2}
  \leq \capacity_{1/K}(F,E_t;B)
  \leq C_3(\diam E'_t)^s.
\end{equation}
Next, by combining~\eqref{eq:epayht1} and~\eqref{eq:valepayht} we obtain that for all $0<t<t_0$
\begin{equation*}
  C_4(\diam E'_t)^{\frac{2+2\varepsilon}{\lambda}}
  \leq C_3(\diam E'_t)^s.
\end{equation*}
From this it follows by taking $t\to0$ (and thus $\diam E'_t\to0$), that for all $\varepsilon>0$ we must have $s\leq(2+2\varepsilon)/\lambda$.But this is a contradiction since $s>0$ is arbitrary. This proves the claim.
\end{proof}
\begin{theorem}\label{thm:uusithm1}
for any $p>0$, there is a homeomorphism $f\colon\bR^2\to\bR^2$  of finite distortion with $K_f(x)\in L^p_{loc}(\bR^2)$ such that $f(B)=\Omega$.
\end{theorem}
\begin{proof}
The construction of the desired mapping $f$ and the computations regarding the distortion are similar to the ones presented in~\cite{kt07, t07}. For the convenience of the reader we represent the entire construction and only omit some analogues computations for the distortion estimate.

The idea behind the construction is that the point $(g(r), e^{-\frac{1}{g(r)}})$ will form the right cusp. So if we express it in polar coordinate, we will get the corresponding mapping to get the right cusp, \ie let $(g(r),e^{-\frac{1}{g(r)}})=(G(r)\cos\theta, G(r)\sin\theta)$ and we obtain
\begin{equation*}
  G(r)= g(r)\sqrt{1+\frac{1}{g(r)}e^{-\frac{1}{g(r)}}},\    \theta = \arctan{H(r)}
\end{equation*}
where $ H(r)= \frac{1}{g(r)}e^{-1/g(r)}$ and $g(r)$ is a test function.

We begin by mapping the unit disc~$B$ conformally onto the open right half plane $H_R:=\{(x_1,x_2)\in\bR^2:x_1>0\}$ so that the point $(-1,0)\in\bdary{B}$ maps to the origin. Using the complex notation we define this mapping $f_1\colon\closure{\bC}\to\closure{\bC}$ as $f_1(z)=(z+1)/(1-z)$.

Next, we define two linear functions $L_r^i\colon(-\frac{\pi}{2},\frac{\pi}{2})\to\bR$ and $L_r^o\colon[\frac{\pi}{2},\frac{3\pi}{2}]\to\bR$ by setting
\begin{align*}
  L_r^i(\theta) &=\frac{2\theta}{\pi}\arctan H(r)\qquad\text{and}\\
  L_r^o(\theta) &=2\theta-\pi+\big(2-\frac{2\theta}{\pi}\big)\arctan H(r).
\end{align*}
Let $x\in\bR^2$ satisfy $|x|\leq1$. We may represent~$x$ in polar coordinates, $x=(r\cos\theta,r\sin\theta)$, so that $0\leq r\leq1$ and either $\theta\in(-\frac{\pi}{2},\frac{\pi}{2})$ or $\theta\in[\frac{\pi}{2},\frac{3\pi}{2}]$. Using this notation, we define the mapping $f_2\colon\bR^2\to\bR^2$ in polar coordinates on~$\closure{B}$ by setting
\begin{equation}\label{eq:f2definB}
f_2(r,\theta)=
\begin{cases}
  (G(r),L_r^i(\theta)) & \text{if $0<r\leq1$ and $\theta\in(-\frac{\pi}{2},\frac{\pi}{2})$} \\
  (G(r),L_r^o(\theta)) & \text{if $0<r\leq1$ and $\theta\in[\frac{\pi}{2},\frac{3\pi}{2}]$} \\
  0 & \text{if $r=0$},
\end{cases}
\end{equation}
where $G(r)=g(r)\sqrt{1+H^2(r)}$.
Outside the closed unit disc~$\closure{B}$, the mapping~$f_2$ will be defined in a \bilip\ manner. First we define $h\colon S(0,1)\to S(0,G(1))$ by setting $h(x)=f_2(x)$ on the unit circle $S:=S(0,1)=\bdary{\closure{B}(0,1)}$. The mapping~$h$ is clearly a \bilip\ mapping on~$S$.
Next we set
\begin{equation}\label{eq:f2defoutB}
f_2(x)=|x|h\big(\frac{x}{|x|}\big) \qquad \text{if } x\in\bR^2\setminus\closure{B}.
\end{equation}
A simple calculation shows that if~$h$ is an $L$-\bilip\ mapping, then $f_2$ will also be an $L$-\bilip\ mapping on $S(0,R)$ for all~$R\geq1$. This and the fact that $|f_2(x)|=|x|G(1)$ for all $|x|\geq1$ assures that~$f_2$ will be a \bilip\ mapping on $\bR^2\setminus\closure{B}$ and the \bilip\ constant of~$f_2$ depends only on~$L$ and~$G(1)$.

The definition in~\eqref{eq:f2definB} gives a mapping that maps the line segment $\{(x_1,x_2)\in\bR^2 : x_1=0,\,-1\leq x_2\leq1\}$ to the set $\{(x_1,x_2)\in\bR^2 : |x_2|=e^{-\frac{1}{x_1}},\,0\leq x_1\leq g(1)\}$, thus forming the desired cusp.

Next, we make this cusp domain bounded by mapping the right half plane onto the disc $B((1/2,0),1/2)$ with the mapping $f_3\colon\closure{\bC}\to\closure{\bC}$, $f_3(z)=z/(z+1)$ and denote $\tilde{\Omega}:=f_3(f_2(f_1(B)))$. The mapping~$f_3$ will somewhat alter the shape of the cusp at the origin, but not essentially as it will be seen. On the other hand, being conformal, it will preserve the size of the two aforementioned non-zero angles and will create one more non-zero angle at the point~$(1,0)$. Otherwise the boundary curve will still be smooth. Because of these facts, there exists a sense preserving \bilip\ mapping~$f_4\colon\bR^2\to\bR^2$ for which $f_4(\tilde{\Omega})=\Omega$ and moreover: this mapping~$f_4$ can be chosen so that for some bounded function $L(x,y)\colon\bR^2\to[1,L]$, $1<L<\infty$, we have that 
\begin{equation}\label{eq:bilipkorj}
  \frac{1}{L(x,y)}|x-y|\leq|f_4(x)-f_4(y)|\leq L(x,y)|x-y|
\end{equation}
and $L(x,y)\to1$ when $|x|+|y|\to0$.
To justify this last claim we notice that by substituting $z=t+ie^{-1/t}$ to $z/(z+1)$, the boundary of cusp takes the form $(t+o(|t|^2),e^{-1/t})$  for small $t$.

Finally we will set $f\colon\bR^2\to\bR^2$, $f(x)=f_4\circ f_3\circ \tilde{f}_2\circ f_1(x)$. Here~$\tilde{f}_2$ is the homeomorphic extension of $f_2$ to~$\closure{\bC}$ obtained by setting~$\tilde{f}_2(x)=f_2(x)$ when~$x$ is finite and $\tilde{f}_2(\infty)=\infty$. This definition clearly gives us a sense preserving homeomorphism for which $f(B)=\Omega$.

Next we show that the distortion function of~$f$ satisfies the required conditions. In fact, it will be enough to compute the distortion of~$f_2$, because the conformal mappings~$f_1$ and~$f_3$ do not give any contribution to it and the contribution of the \bilip\ mapping~$f_4$ will be in fact a bounded multiplier that goes to one when we approach the origin. The latter claim follows from the facts that the distortion of an $L$-\bilip\ mapping is $L^2$, and~$f_4$ satisfies the inequality~\eqref{eq:bilipkorj}.

Outside the unit disc~$B$ the mapping~$f_2$ is quasiconformal, as it is seen directly from the definition to be sense preserving and \bilip\ there. 
Hence we need to show that~$f_2$ has an $L^p$-integrable distortion in~$B$.

As computed in~\cite{kt07, t07}, for the case $\theta\in [0, \pi/2]$, the resulting differential matrix is
\begin{multline}\label{eq:dmatriisi}
\begin{bmatrix}
 \frac{d}{dr}G_s(r) & 0 \\
 G_s(r)\frac{d}{dr}L_r^i(\theta) & \frac{G_s(r)}{r}\frac{d}{d\theta}L_r^i(\theta)\rule{0pt}{20pt}
\end{bmatrix}\\
=\begin{bmatrix}
  \scriptstyle\frac{1}{\sqrt{1+H^2(r)}}[g'(r)(1+H^2(r))+g(r)H(r)H'(r)] & \scriptstyle 0 \\
  \scriptstyle\frac{2\theta}{\pi}\frac{1}{\sqrt{1+H^2(r)}}g(r)H'(r) & %
  \scriptstyle\frac{2}{\pi}\frac{g(r)}{r}\sqrt{1+H^2(r)}\arctan H(r)\rule{0pt}{14pt}
 \end{bmatrix};
\end{multline}
for the case $\theta\in[\frac{\pi}{2},\frac{3\pi}{2}]$, the differential matrix
\begin{multline*}
\begin{bmatrix}
 \frac{d}{dr}G_s(r) & 0 \\
 G_s(r)\frac{d}{dr}L_r^o(\theta) & \frac{G_s(r)}{r}\frac{d}{d\theta}L_r^o(\theta)\rule{0pt}{20pt}
\end{bmatrix}\\
=\begin{bmatrix}
  \scriptstyle\frac{1}{\sqrt{1+H^2(r)}}[g'(r)(1+H^2(r))+g(r)H(r)H'(r)] & \scriptstyle 0 \\
  \scriptstyle\frac{2-2\theta/\pi}{\sqrt{1+H^2(r)}}g(r)H'(r) & %
  \scriptstyle(2-\frac{2}{\pi}\arctan H(r))\frac{g(r)}{r}\sqrt{1+H^2(r)}\rule{0pt}{14pt}
 \end{bmatrix}.
\end{multline*}

Now, if we choose the test function $g(r)=\frac{1}{\log{\log{\frac{2}{r}}}}$, then a direct calculation shows
\begin{equation}\label{eq:arvioulkona}
  K_{f_2}(x)\leq C \log \frac{2}{|x|}\cdot\log{\log{\frac{2}{|x|}}}
\end{equation}
The estimate ~\eqref{eq:arvioulkona} show that $K_{f_2}(x)\in L^p(B)$ for all $p>0$.

As mentioned before, $f_1$ and~$f_3$ do not give any contribution to the distortion of~$f$ and the contribution of~$f_4$ is a bounded multiplier that goes to one when one approaches the origin. Moreover, $f_1$ is \bilip\ in~$f^{-1}(B)$. Thus the distortion~$K_f(x)$ of~$f$ satisfies $K_f(x)\in L^p_{\loc}(\bR^2)$ for all $p>0$.
\end{proof}
\paragraph{\textbf{Acknowledgements.}} I would like to thank my advisor Professor Pekka Koskela for his guidance and comments.

\end{document}